\documentclass[letterpaper, 10 pt, conference]{ieeeconf}   \IEEEoverridecommandlockouts                              
\usepackage{amssymb}
\usepackage{amsmath}
\usepackage{graphicx}
\usepackage{xspace}
\usepackage[FIGTOPCAP]{subfigure}
\let\labelindent\relax
\usepackage{enumitem}   
\usepackage{multirow}
\usepackage{cite}
\usepackage{amsfonts}
\usepackage{mathrsfs}

\usepackage{url}

\usepackage{theorem}

\newtheorem{theorem}{Theorem}
\newtheorem{proposition}{Proposition}

\newtheorem{definition}{Definition}
\newtheorem{remark}{Remark}

\usepackage[final]{listofsymbols}
\usepackage{algorithm}
\usepackage[ruled,vlined,algo2e,linesnumbered]{algorithm2e}
\usepackage[overload]{empheq}

\usepackage[table]{xcolor}      % Provides coloring for tables and text

\usepackage[acronym]{glossaries}

\allowdisplaybreaks[4]
\newacronym{ocd}{OCD}{Optimality Condition Decomposition}
\newacronym{app}{APP}{Auxiliary Problem Principle}
\newacronym{admm}{ADMM}{Alternating Direction Method of Multipliers}
\newacronym{aladin}{ALADIN}{Augmented Lagrangian based Alternating Direction Inexact Newton method}
\newacronym{rapidpf}{\texttt{rapidPF}}{rapid prototyping for distributed power flow}
\newacronym{rapidpf+}{\texttt{rapidPF}+}{rapid prototyping for distributed power flow}
\newacronym{bfgs}{BFGS}{Broyden-Fletcher-Goldfarb-Shanno}
\newacronym{tsos}{TSOs}{transmission system operators}
\newacronym{dsos}{DSOs}{distribution system operators}
\newacronym{qp}{QP}{quadratic programing}
\newacronym{sqp}{SQP}{sequential quadratic programing}
\newacronym{dsqp}{DSQP}{distributed sequential quadratic programing}
\newacronym{hdqp}{HDQ}{hypergraph-based distributed quadratic optimization algorithm}
\newacronym{qcqp}{QCQP}{Quadratically Constrained Quadratic Program} 
\newacronym{licq}{LICQ}{linear independence constraint qualification} 
\newacronym{sosc}{SOSC}{second order sufficient condition} 
\newacronym{scc}{SCC}{strict complementarity conditions} 
\newacronym{kkt}{KKT}{Karush–Kuhn–Tucker} 
\newacronym{pf}{PF}{power flow}
\newacronym{opf}{OPF}{optimal power flow}
\newacronym{hdsqp}{HDSQP}{hypergraph-based distributed sequential quadratic programming}
\newacronym{nlp}{NLP}{Nonlinear Programming}
\newcommand{\change}{\textcolor{black}}
\newcommand{\matlab}{\textsc{matlab}\xspace}
\newcommand{\norm}[1]{\left\lVert#1\right\rVert}

\newcommand{\matpower}{\textsc{matpower}\xspace}

\newcommand{\ipopt}{\textsc{ipopt}\xspace}

\newcommand{\casadi}{\textsc{c}as\textsc{ad}i\xspace}

\newcommand{\ac}{\textsc{ac}\xspace}
\def\subparagraph{} % because IEEE classes don't define this, but titlesec assumes it's present
\usepackage{titlesec}
\titlespacing*{\section}{0pt}{*0.8}{*0.8}
\titlespacing{\subsection}{0pt}{*0.8}{*0.8}
\renewcommand{\thesubsubsection}{\arabic{subsubsection}}
\titleformat{\subsubsection}[runin]{\itshape}{\thesubsubsection)}{1em}{}
\titlespacing*{\subsubsection}{\parindent}{0pt}{*0.8}

\title{\LARGE \bf
Hypergraph-Based Fast Distributed AC Power Flow Optimization%*
}

\author{Xinliang Dai, Yingzhao Lian, Yuning Jiang, Colin N. Jones and Veit Hagenmeyer
\thanks{ This work was supported in part by the German Federal Ministry of Education and Research within the project MOReNet – Modellierung, Optimierung und Regelung von Netzwerken heterogener Energiesysteme mit volatiler erneuerbarer Energieerzeugung, in part by the BMBF-project ENSURE II with grant number 03SFK1F0-2 and in part by the Swiss National Science Foundation (SNSF) under the NCCR Automation project, grant agreement 51NF40\_180545. }
\thanks{XD and VH are with Karlsruhe Institute of Technology, Karlsruhe, Germany. Email: {\tt\{xinliang.dai, veit.hagenmeyer\}@kit.edu}}
\thanks{YZ, YJ and CJ are with École Polytechnique Fédérale de Lausanne, Lausanne, Switchland. Email: {\tt \{yingzhao.lian, yuning.jiang, colin.jones\}@epfl.ch}}
}

\begin{document}

\maketitle
\setlength\abovedisplayskip{2.pt}
\setlength\belowdisplayskip{2.pt}

%\thispagestyle{empty}
%\pagestyle{empty}

%%%%%%%%%%%%%%%%%%%%%%%%%%%%%%%%%%%%%%%%%%%%%%%%%%%%%%%%%%%%%%%%%%%%%%%%%%%%%%%%
\begin{abstract}
This paper presents a novel distributed approach for solving AC power flow (PF) problems. The optimization problem is reformulated into a distributed form using a communication structure corresponding to a hypergraph, \change{by which complex relationships between subgrids can be expressed as hyperedges}. Then, a hypergraph-based distributed sequential quadratic programming (\acrshort{hdsqp}) approach is proposed to handle the reformulated problems, \change{and the \acrfull{hdqp} is used as the inner algorithm to solve the corresponding QP subproblems, which are respectively condensed using Schur complements with respect to coupling variables defined by hyperedges. Furthermore, we rigorously establish the convergence guarantee of the proposed algorithm with a locally quadratic rate and the one-step convergence of the inner algorithm when using the Levenberg-Marquardt regularization. Our analysis also demonstrates that the computational complexity of the proposed algorithm is much lower than the state-of-art distributed algorithm.} We implement the proposed algorithm in an open-source toolbox, \texttt{rapidPF}\footnote{Open-source toolkit: \url{https://github.com/xinliang-dai/rapidPF}}, and conduct numerical tests that validate the proof and demonstrate the great potential of the proposed distributed algorithm in terms of \change{communication effort} and computational speed.
\end{abstract}

\section{Introduction}
\label{sec::intro}
The \acrfull{pf} problem is a fundamental problem in power system analysis and has many applications, such as planning, expanding, and operating power systems~\cite{grainger1999power}. \change{Traditionally, centralized methods such as Gauss-Seidel~\cite{glimm1957seidel} or Newton-type methods~\cite{tinney1967newtonpf,stott1974newtonpf} have been used to solve PF problems. In recent years, several studies were carried out in various aspects, including analysis of power flow equations~\cite{dvijotham2015powerflow}, state estimation~\cite{guo2023online,guo2023optimal}, distributionally robust optimal
control~\cite{guo2016case}, initialization strategies~\cite{thorp1989initial,thorp1990initial}, convex relaxation~\cite{madani2015convexification,molzahn2019survey}, and convex restriction~\cite{lee2019restriction}.} With the increasing penetration of distributed energy resources and the need for optimization and control of power systems with many controllable devices, distributed approaches have gained significant research attention~\cite{molzahn2017survey}. For systems like Germany's power grid, which has four \acrfull{tsos} and over 900 \acrfull{dsos}, sharing detailed grid models is not preferred. Therefore, a centralized approach is not preferred by systems operators or is even prohibited by the respective regulation.

The present paper focuses on AC models to obtain more realistic results. The main challenge is that the AC \acrshort{pf} feasibility is NP-hard~\cite{lehmann2015ac,bienstock2019strong}, and is a challenge even for a centralized approach. \cite{sun2008distributed} proposed to solve AC \acrshort{pf} problems by breaking the original problem into several smaller power flow subproblems, keeping coupling variables fixed, and then iterating over them. In the follow-up work~\cite{sun2014master}, the convergence was analyzed under some additional assumptions. However, the actual convergence behaviors and scalability are limited in practice. Other well-known distributed algorithms for AC power flow analysis---\acrfull{ocd} proposed by~\cite{hug2009decentralized}, \acrfull{app} by~\cite{baldick1999fast}, and \acrfull{admm} by~\cite{erseghe2014distributed}---have no convergence guarantees in general, and their convergence behaviors are case-by-case in practice.

Recently, \cite{houska2016augmented}  proposed a second-order distributed algorithm, i.e., \acrfull{aladin}. In contrast to these existing distributed approaches, \acrshort{aladin} can provide a local convergence guarantee with a quadratic convergence rate for generic distributed nonconvex optimization problems if suitable Hessian approximations are used. Based on the \acrshort{aladin} algorithm, considerable works have been carried out for power system analysis~\cite{engelmann2018toward, ZhaiJunyialadin,bauer2022shapley}. \cite{muhlpfordt2021distributed} provides open-source \matlab code for \acrfull{rapidpf}. Extensive research~\cite{dai2022rapid} has improved computing time significantly for solving large-scale AC \acrshort{pf} problems by using Gauss-Newton approximation and further exploiting the problem formulation. However, \acrshort{aladin} is limited by the required computation and communication effort per iteration.

The aforementioned studies either intertwine problem formulation and problem solution or use the standard affinely coupled distributed form. In contrast, we introduce a hypergraph~\cite{berge1973graphs} based AC \acrshort{pf} framework in the present paper. \change{As a generalization of graphs, hypergraph allows more than two nodes to be connected in the same hyperedge, therefore depicting more complex relationships}, e.g., multiple regions connected to a region at the same bus. More specifically, we propose to reformulate the AC \acrshort{pf} problem as a zero-residual least-squares problem with a communication structure corresponding to a hypergraph and then solve it by an \acrfull{hdsqp} approach. The convergence of the proposed \acrshort{hdsqp} is guaranteed, and the convergence rate is quadratic when approaching the minimizer if the Levenberg-Marquardt method is used. Moreover, the condensed \acrshort{qp} subproblems of \acrshort{hdsqp} are solved by a variant of the dual decomposition algorithm, i.e.,  the \acrfull{hdqp} proposed by~\cite{papastaikoudis_hypergraph_2022}. Most notably, the communication matrix of the dual decomposition is the Bollas' Laplacian for hypergraphs~\cite{bolla1993spectra}, and it is also a projection matrix for the \acrshort{hdqp} algorithm. As a result, \acrshort{hdqp} can converge to the global minimizer in one iteration if the corresponding \acrshort{qp} subproblem is convex. This hints at the fact that the proposed distributed approach could converge rapidly with less communication effort.

%The aim of the present paper is to exploit the hypergraph-based distributed approach for solving generic AC \acrshort{pf} problems. To this end,  Section~\ref{sec::formulation} proposes formulating the AC \acrshort{pf} problem as a zero-residual least-squares problem with a communication structure corresponding to a hypergraph. The main contribution of this work, presented in Section~\ref{sec::algorithm}, is a novel \acrfull{hdsqp} approach for solving the  reformulated and structured AC \acrshort{pf} problems proposed in Section~\ref{sec::formulation}. In particular, we tailor the classical \acrfull{sqp} framework, and the \acrfull{hdqp} is used as the inner algorithm to solve the corresponding \acrshort{qp} subproblems. Furthermore, Section~\ref{sec::analysis} rigorously analyzes the local convergence of the proposed algorithm. In particular, Theorem~\ref{thm::gn::rate} establishes the local quadratic convergence rate of the proposed \acrshort{hdsqp} if the Levenberg-Marquardt method is used. Numerical tests presented in  Section~\ref{sec::simulation} show that the proposed algorithm surpasses the state-of-art \acrshort{aladin} algorithm with respect to computing time and communication effort. Moreover, our algorithm has been implemented as one module in the open-source Toolbox \texttt{rapidpf}. Finally, we conclude the paper in Section~\ref{sec::conclusion}.

The aim of the present paper is to exploit the hypergraph-based distributed approach for solving generic AC \acrshort{pf} problems. The main contributions are listed in the following:
\begin{enumerate}[label=(\roman*)]
     \item We propose a new distributed form of the AC \acrshort{pf} problem that uses a communication structure corresponding to a hypergraph to generalize complex relationships between subgrids. Moreover, we propose a \acrfull{hdsqp} approach for solving the problem and use the \acrfull{hdqp}~\cite{papastaikoudis_hypergraph_2022} as the inner algorithm for the corresponding condensed \acrshort{qp} subproblem at each iteration.      
     \item We rigorously establish the convergence guarantee of the proposed algorithm \acrshort{hdsqp} with a locally quadratic rate and the one-step convergence of the inner algorithm \acrshort{hdqp} when using the Levenberg-Marquardt regularization. Our analysis also demonstrates that the computational complexity of the proposed algorithm is much lower than the state-of-art distributed algorithm. Numerical tests are added to the \acrshort{rapidpf} open-source toolbox, and we
     show that the proposed algorithm surpasses the state-of-art \acrshort{aladin} algorithm with respect to computing time and communication effort. 
 \end{enumerate} 
 
This paper introduces the distributed formulation of AC \acrshort{pf} problem in Section~\ref{sec::formulation}. Then, we present the proposed hypergraph-based distributed optimization algorithm and the convergence analysis in Section~\ref{sec::algorithm}. Finally, we present numerical simulations in Section~\ref{sec::simulation} and conclude the paper in Section~\ref{sec::conclusion}. Additionally, an anonymous chat is listed in Appendix.

\section{System Model and Problem Formulation}
\label{sec::formulation}
This paper considers a power system defined by a tuple $\mathcal{S}=(\mathcal{R},\;\mathcal{N},\;\mathcal{L})$ with the set $\mathcal{R}$ of all regions, $\mathcal{N}$ the set of all buses and $\mathcal{L}$ the set of all branches. We define by $ n^\textrm{reg}$, $n^\textrm{bus}$, and $n^\textrm{line}$ the number of regions, buses, and branches, respectively.  In the present paper, we use complex voltage in polar coordinates: 
$$
V_i= v_i e^{\textrm{j}\theta_i},\;\;i\in\mathcal{N}.
$$ 
where $v_i$ and $\theta_i$ denote the voltage magnitude and angle. Thereby, for each bus $i\in\mathcal N$, its steady state $(\theta_i,v_i,p_i,q_i)$
includes $v_i$ and $\theta_i$ the voltage magnitude and angle, $p_i$ and $q_i$ the active and reactive power. Throughout this paper, we stack all steady states at region $\ell\in\mathcal R$ by $\chi_\ell=\{(\theta_i,v_i,p_i,q_i)\}_{i\in\mathcal N_\ell}$ with $\mathcal N_\ell$ bus set at region $\ell$. When a vector $\xi$ consists of $ n^\textrm{reg}$ subvectors, we write $\xi =(\xi_1,...,\xi_{ n^\textrm{reg}})$. Moreover, $x_\ell$ denotes the coupled variables, and $y_\ell$ denotes hidden variables that are totally local to region $\ell\in\mathcal{R}$, i.e., $\chi_\ell := (x_\ell,y_\ell)$. Thus, there is a matrix $A_\ell$ such that 
$x_\ell = A_\ell  \chi_\ell$.
\subsection{Hyergraph-based modeling}%Moreover, we use $B_\ell$ to represent the Hessian approximation of a specific subproblem $\ell\in\mathcal{R}$, and $B$ to represent the corresponding block diagonal Hessian approximation for the coupled problems.

\begin{figure}[htbp!]
    \begin{center}
    \subfigure[Multi-region coupled system]{
        \includegraphics[width=0.25\textwidth]{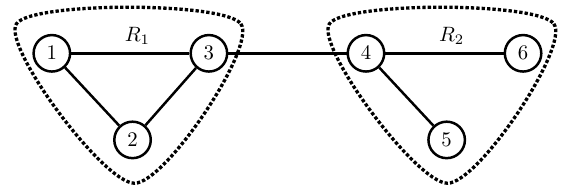}
        \label{fig::coupled}
    }\\
    \subfigure[Hypergraph]{
        \includegraphics[width=0.25\textwidth]{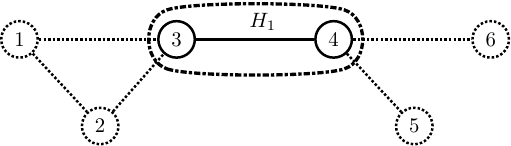}
        \label{fig::hypergraph}
    }\\
    \subfigure[Decoupled region 1]{
        \includegraphics[width=0.2\textwidth]{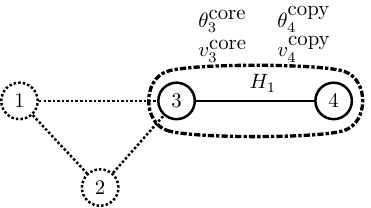}
        \label{fig::region1}
    }
    \hskip 20pt
    \subfigure[Decoupled region 2]{
        \includegraphics[width=0.2\textwidth]{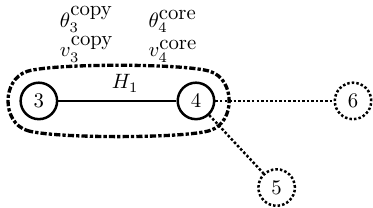}
        \label{fig::region2}
    }
    \end{center}
    \caption[Optional caption for list of figures]{Decomposition by sharing components for a two-region system}\label{fig::example}
\end{figure}

As discussed in~\cite{muhlpfordt2021distributed}, we share the components between neighboring regions to ensure physical consistency. Let us take the 6-bus system with two regions, shown in Fig.~\ref{fig::example} as an example. The coupled system, shown in Fig.~\ref{fig::example}(a), has been partitioned into two local regions. To solve the \ac \acrshort{pf} problem in the region $R_1$, besides its buses \{1,2,3\} called the \textit{core buses}, the complex voltage of bus \{4\} from neighboring region $R_2$ is required. Hence, for the sub-problem of the region $R_1$, we create an auxiliary bus \{4\} called the \textit{copy bus}, along with its own \textit{core bus}, to formulate a self-contained AC \acrshort{pf} problem. The resulting affine consensus constraint can be written as
\begin{subequations}\label{eq::consensus::original}
    \begin{align}
        \theta^{\textrm{core}}_3 = \theta^{\textrm{copy}}_3,&\;\theta^{\textrm{core}}_4 = \theta^{\textrm{copy}}_4,\\
        v^{\textrm{core}}_3 = v^{\textrm{copy}}_3,&\;v^{\textrm{core}}_4 = v^{\textrm{copy}}_4.
    \end{align} 
\end{subequations}
The multi-region power system $\mathcal{S}$ can be transformed into a hypergraph $\mathcal{G}=(\mathcal{N}, \mathcal{H})$, where $\mathcal{H}=(H_1,\cdots,H_n)$ denotes the set of all hyperedges that can include any number of buses, c.f. \cite{voloshin2009introduction}. The example mentioned above has only one hyperedge, as shown in Fig.~\ref{fig::example}(b). The buses $\{1,2,5,6\}$ are called isolated or hidden buses since they are not incident to any hyperedge. Thereby, we have coupling variables of region $R_1$ and $R_2$, i.e., $x_1 = (\theta_3^\textrm{core}, \theta_4^\textrm{copy},v_3^\textrm{core},v_4^\textrm{copy})$ and  $x_2 = (\theta_3^\textrm{copy}, \theta_4^\textrm{core},v_3^\textrm{copy},v_4^\textrm{core})$. Moreover, we denote by $z =(\theta_3, \theta_4,v_3,v_4)$ to represent the common values for coupling variables of all regions. As a result, the consensus constraints~\eqref{eq::consensus::original} can be written as
\begin{equation}
    x_\ell = E_\ell z,\;\ell\in\mathcal{R},
\end{equation}
where $E_\ell$ denote the incident matrix of a specific region $\ell\in\mathcal{R}$. 
\subsection{Hypergraph-based AC Power Flow Problem}
In polar coordinates, the resulting conventional AC \acrfull{pf} problem can be written as a set of power flow equations
\begin{subequations}\label{eq::pf::equation}
    \begin{align}
    p_i = p_i^g-p_d^l &=v_i \sum_{k\in\mathcal{N}} v_k \left( \mathcal G_{ik} \cos\theta_{ik} +\mathcal B_{ik} \sin\theta_{ik} \right),\\
    q_i = q_i^g-q_d^l &=v_i \sum_{k\in\mathcal{N}} v_k \left( \mathcal G_{ik} \sin\theta_{ik} -\mathcal B_{ij} \cos\theta_{ik} \right)
    \end{align}
\end{subequations}
for all buses $i \in\mathcal{N}$, where $p^g_i$, $q^g_i$ (reps. $p^d_i$, $q^d_i$) denote the real and reactive power injections from generator(s) (resp. loads) at bus $i$, $\theta_{ik}$ denotes the angle difference between bus $i$ and bus $k$, $\mathcal G_{ij}$, $\mathcal B_{ij}$ are the real and imaginary parts of the bus admittance matrix entries $Y_{ik} =\mathcal G_{ij} +  \textrm{j} \mathcal B_{ij}$. 
\begin{remark}
    Since multiple stable AC power flow solutions may exisit~\cite{korsak1972question,mehta2016recent}, especially in the presence of power flow reversal on distribution systems~\cite{nguyen2014appearance,nguyen2015voltage}, the present paper focuses on local solutions.
\end{remark}
Following~\cite{dai2022rapid}, these power flow equations can be written as a residual function
\begin{equation}\label{eq::pf::residual}
    r(\chi)=0
\end{equation}
with state variables $\chi = (\theta,v,p,g)$, so that the AC \acrshort{pf} problems can be formulated as a standard zero-residual least-squares problem
\begin{equation}\label{eq::pf::ls::original}
    \min_\chi\; f(\chi) = \frac{1}{2} \norm{r(\chi)}_2^2.
\end{equation}
In other words, $f$ is the sum of squared residuals of power flow equations for all buses $i\in\mathcal{N}$. Hence, both the state variables $\chi$ and the objective $f$ are separable, i.e.,
\begin{equation}
    f(\chi) = \sum_{\ell\in\mathcal{R}} f_\ell(\chi_\ell)=\sum_{\ell\in\mathcal{R}} \frac{1}{2}\|r_\ell(\chi_\ell)\|_2^2.
\end{equation}
As a result, the coupled problem~\eqref{eq::pf::ls::original} can be reformulated with a communication structure corresponding to a hypergraph
\begin{subequations}\label{eq::pf::ls::hypergraph}
    \begin{align}
        \min_\chi\quad & f(\chi)= \sum_{\ell\in\mathcal{R}} f_\ell(\chi_\ell)\\\label{eq::coupling}
        \text{subject to}  \quad  & x_\ell = \; E_\ell \,z\quad\mid\lambda_\ell,\;\; \ell\in\mathcal{R}.
    \end{align}
\end{subequations}
Recall that $\chi$ consists of two components, i.e., 
\begin{equation}
    \chi_\ell = (x_\ell,\;y_\ell),\quad \forall \ell\in\mathcal{R},
\end{equation}
where $x_\ell$ represents coupling variables and $y_\ell$ represents hidden variables that are entirely local.

\section{Distributed Optimization Algorithm}\label{sec::algorithm}
This section introduces the \acrfull{sqp} framework and the condensed reformulation of \acrshort{qp} subproblem at each iteration. Then, based on the preliminaries introduced in Section~\ref{sec::formulation}, we propose a hypergraph-based distributed approach to solve the AC \acrshort{pf} problem~\eqref{eq::pf::ls::hypergraph}. Convergence analysis is carried out at the end of this section.

\subsection{Preliminaries}
\label{sec::pre}
The sequential quadratic programming~(\acrshort{sqp}) framework is one of the most effective methods for \acrfull{nlp}, in which a sequence of \acrshort{qp} subproblems are iteratively constructed and solved. This paper focuses on the full-step variant and its corresponding local convergence. Because in practice, a conventional flat start can always provide a good initial guess for AC \acrshort{pf} problems. Regarding the globalization routine or line search for \acrshort{sqp}, more details refer to~\cite[Chapter~18]{nocedal2006numerical}. 

At the $k$-th iteration, the \acrshort{qp} subproblem is written as 
\begin{subequations}\label{eq::subproblem}
    \begin{align}
        \min_{\chi^{k+1},z}\; &\;m^k(\chi^{k+1})=\sum_{\ell\in\mathcal{R}}\, m_\ell^k(\chi_\ell^{k+1})\\
        \text{subject to}\;&\; A_\ell \chi_\ell^{k+1}= E_\ell z\quad|\quad\lambda_\ell,\;\forall\ell\in\mathcal{R}.
    \end{align}  
\end{subequations}
with quadratic models
\begin{equation}\label{eq::qp::model}
    m^k_\ell(\chi_\ell) = \frac{1}{2} (\chi_\ell)^\top \nabla^2 f_\ell^k\, \chi_\ell + \left(\nabla f^k_\ell - \nabla^2 f^k_\ell \,\chi^k_\ell\right)^\top \chi_\ell,
\end{equation}
and for notation simplication, $f^k_\ell = f_\ell(\chi^k_\ell)$ used for all $\ell\in\mathcal{R}$. 
The derivatives of the objectives $f_\ell$ at iterate $\chi^k_\ell$ can be expressed as
\begin{subequations}
    \begin{align}
        \nabla f_\ell^k&= \left(J_{\ell}^k\right)^{\top} r_\ell^k,\\
        \nabla^2 f_\ell^k &= \left(J_{\ell}^k\right)^{\top} J_{\ell}^k +  Q_{\ell}^k\label{eq::second::derivative}
    \end{align}
\end{subequations}
with 
\begin{subequations}
    \begin{align}
    J_\ell^k &= J_\ell(\chi_\ell^k) = 
    \begin{bmatrix}
    \nabla r_{_\ell,1},\nabla r_{_\ell,2},\cdots,\nabla r_{_\ell,n_\textrm{pf}}
    \end{bmatrix}^{\top},\label{eq::hessian::first}\\
    Q_\ell^k &= Q_{\ell}(\chi_{\ell}^k) = \sum_{m=1}^{n_\textrm{pf}}r_{\ell,m}(\chi_{\ell}^k) \nabla^2 r_{\ell,m}(\chi_{\ell}^k).\label{eq::hessian::second}
    \end{align}
\end{subequations}
Here, $r_{\ell,m}$ represents the residual of the $m$-th power flow equation at region $\ell\in\mathcal{R}$. In practice, the first term $\left(J^k_\ell\right)^\top\;J^k_\ell$ dominates the second term $Q^k_\ell$, because the residuals $r_{\ell,m}$ are close to zero near the solution~\cite[Chapter~10]{nocedal2006numerical}. In the present paper, Levenberg-Marquardt regularization is used, i.e.,
\[
B^k_\ell = \left(J^k_\ell\right)^\top J^k_\ell+ \varepsilon I,\quad\forall\ell\in\mathcal{R},
\]
to approximate Hessians such that~\eqref{eq::subproblem} is strongly convex. Here, one empirical choice of $\varepsilon$ in practice, is $\varepsilon = 10^{-10}$. As discussed in~\cite[Chapter~10]{nocedal2006numerical}, when $\varepsilon$ is sufficiently small, the Levenberg-Marquardt method shares the same performance with the classical Gauss-Newton method under mild assumptions. 
%\begin{rem}\label{rem::inexact}
%For solving the zero-residual least-squares problem, the Gauss-Newton approximation is exact at the solution, and can converge to the exact Hessian during iterations rapidly.
%\end{rem}
%
%The Gauss-Newton  is symmetric positive semidefinite.  %Moreover, derivatives of the condensed \acrshort{qp} subproblems~\eqref{eq::dsqp::subqp::reduced}, i.e., $\overline{B}^k_\ell\approx\nabla^2 F^k_\ell$ and $\overline{g}^k_\ell = \nabla F^k_\ell$, can be evaluated using the Schur complement~\eqref{eq::schur} by replacing $\nabla^2 f^k_\ell$ by Levenberg-Marquardt approximation $B^k_\ell$.

In order to write the \acrshort{qp} subproblems~\eqref{eq::subproblem} in a condensed form, we set
\begin{equation}
    F_\ell(x_\ell) = \min_{y_\ell} f_\ell(x_\ell,y_\ell),
\end{equation}
where $x_\ell$ denotes the coupling variables and $y_\ell$ denotes the variables that are entirely local. Accordingly, the Jacobian matrix $J^k_\ell$ can be partitioned into two blocks w.r.t. $x_\ell$ and $y_\ell$ as
\begin{equation}
    J^k_\ell = \begin{bmatrix}
    J^x_\ell,\;J^y_\ell
    \end{bmatrix},\quad \forall \ell\in\mathcal{R}.
\end{equation}
Consequently, the gradient $g^k_\ell$ and the approximated Hessians $B^k_\ell$ can be written as
\begin{equation}
    \hspace{-2mm} 
    g^k_\ell = \begin{bmatrix}
    \left(J^x_\ell\right)^\top r^x_\ell,\;\left(J^y_\ell\right)^\top r^y_\ell
    \end{bmatrix}\;\textrm{and}\;B^k_\ell = \begin{bmatrix}
    B^{xx}_\ell & B^{xy}_\ell \\[0.12cm]
    B^{yx}_\ell & B^{yy}_\ell
    \end{bmatrix}
\end{equation}
with the assistance of the Levenberg-Marquardt method
\begin{subequations}
\begin{align}
    B^{xx}_\ell &= \left(J^x_\ell\right)^\top J^x_\ell + \varepsilon I,\\
    B^{yy}_\ell &= \left(J^y_\ell\right)^\top J^y_\ell + \varepsilon I,\\
    B^{xy}_\ell &= \left(B^{xy}_\ell\right)^\top =  \left(J^x_\ell\right)^\top J^y_\ell
\end{align}
\end{subequations}
for all $\ell\in\mathcal{R}$. By using the Schur complement, the first and the second derivatives of the function, $F_\ell$ can be given by
\begin{subequations}\label{eq::schur}
    \begin{align}
        \overline{g}^k_\ell&= g^x_\ell - B^{xy}_\ell \left[B^{yy}_\ell\right]^{-1} g^y_\ell,\\[0.12cm]
        \overline{B}^k_\ell &= B^{xx}_\ell - B^{xy}_{\ell}  \left[B^{yy}_\ell\right]^{-1} B^{yx}_\ell.
        %\overline{c}^k &= \overline{g}^k_\ell - \overline{B}^k_\ell x^k_\ell
    \end{align}
\end{subequations}
Thereby, the \acrshort{qp} subproblems~\eqref{eq::subproblem} can be reformulated into a condensed form
\begin{subequations}\label{eq::dsqp::subqp::reduced}
    \begin{align}
        \min_{x,z}\quad & \overline{m}^k(x)= \sum_{\ell\in\mathcal{R}} \overline{m}^k_\ell(x_\ell)\\
        \text{subject to}\quad & x_\ell = E_\ell\, z\qquad|
        \;\lambda_\ell,\;\; \forall \ell \in \mathcal{R}
    \end{align}    
\end{subequations}
with reduced quadratic models
\begin{equation}\label{eq::reduced::model}
    \overline{m}^k_\ell (x_\ell) = \frac{1}{2} x_\ell^\top \;\overline{B}^k_\ell \;x_\ell + \left(\overline{g}^k_\ell - \overline{B}^k_\ell x^k_\ell\right)^\top x_\ell.
\end{equation}
Here, $x^k_\ell = A_\ell\chi^k $ for all $\ell\in\mathcal{R}$. The Lagrangian function of the condensed \acrshort{qp} subproblem~\eqref{eq::dsqp::subqp::reduced} can thus be written as
\begin{equation}\label{eq::subproblem::lagrangian}
        \mathcal{L}(x,z,\lambda) = \sum_{\ell\in\mathcal{R}} \left\{\overline{m}^k_\ell(x_\ell) + \lambda^\top_\ell x_\ell \right\} - \lambda^\top E\,z
\end{equation}
with $E=[E_1^\top,\cdots,E_\ell^\top]^\top$. The \acrfull{kkt} conditions of~\eqref{eq::dsqp::subqp::reduced} are as follows,
\begin{subequations}\label{eq::kkt}
\begin{align}
    \nabla_x \mathcal{L}  =&\;0\;=\; \overline{B}^k (x -x^k) +  \overline{g}^k + \lambda, \label{eq::kkt::subqp::1}\\        \nabla_z \mathcal{L}  =&\;0\;=\; E^\top \lambda, \label{eq::kkt::subqp::dual}\\
    \nabla_\lambda \mathcal{L}  =&\;0\;=\; x - E\;z, \label{eq::kkt::subqp::primal}
\end{align}
\end{subequations}
where $\overline{B}^k=\text{diag}\{\overline{B}_\ell^k\}_{\ell\in\mathcal R}$ stacks all $\overline{B}_\ell$ into a block diagonal matrix.
\subsection{Hypergraph-Based SQP}

Based on Section~\ref{sec::pre}, we propose an \acrfull{hdsqp} approach to solve~\eqref{eq::pf::ls::hypergraph}. More specifically, the \acrfull{hdqp}, proposed as a variant of dual decomposition by~\cite{papastaikoudis_hypergraph_2022}, is implemented to solve its condensed \acrshort{qp} subproblems~\eqref{eq::dsqp::subqp::reduced}. Remarkably, \acrshort{hdqp} can converge with convexity assumption in one iteration to save total computing time and communication effort. 

Algorithm~\ref{alg} outlines the proposed \acrshort{hdsqp}. Line~\ref{alg::evaluation} evaluates derivatives of the full-dimentional subproblem~\eqref{eq::subproblem} and the corresponding condensed subproblems~\eqref{eq::dsqp::subqp::reduced} with assistance of the Levenberg-Marquardt method~\eqref{eq::lm} and the Schur complement~\eqref{eq::schur}. The resulting condensed subproblem~\eqref{eq::dsqp::subqp::reduced} is a strongly convex \acrshort{qp} with a communication structure corresponding to a hypergraph. Thereby, \acrshort{hdqp} is added as inner algorithm to solve the condensed \acrshort{qp} subproblem~\eqref{eq::dsqp::subqp::reduced} (Line~\ref{alg::hdqp::decoupled}-\ref{alg::hdqp::update}) due to fast convergence rate. The \acrshort{hdqp} algorithm consists of three steps. In Line~\ref{alg::hdqp::decoupled}, temporary local coupling variables $\bar x_\ell$ for all region $\ell\in\mathcal{R}$ are obtained with respect to the KKT condition~\eqref{eq::kkt::subqp::1} under initial condition $\lambda_\ell=0$. Then, weighted averaging is conducted to compute a temporary state $\bar z$ in Line~\ref{alg::hdqp::averaging}, where the weights are determined by the condensed Hessian approximation $\overline{B}^k=\text{diag}\{\overline B_{\ell}^k\}_{\ell\in\mathcal R}$. In the fourth step, the dual variable $\lambda$ is updated based on the deviation of temporary weighted primal residual $\overline{B}^k(\bar x-E\bar z)$ in~\eqref{eq::hdq::dual}. 
\begin{algorithm}[htbp!]
\SetAlgoLined
\textbf{Initialization:} $\chi^0$ as a flat start\\
\Repeat{Primal variables $\chi$ converge}{
    \nl Evaluate derivatives at iterate $\chi^k$\label{alg::evaluation}
    \begin{equation}\label{eq::lm}
        B^k_\ell = \left(J^k_\ell\right)^\top J^k_\ell+\varepsilon I \quad\textrm{and}\quad g^k_\ell = \left(J^k_\ell\right)^\top r^k_\ell,
    \end{equation}\\
    and the corresponding condensed derivatives $\overline{B}_\ell^k$ and $\overline{g}^k_\ell$ by Schur complement~\eqref{eq::schur} for all $\ell\in\mathcal{R}$ \\[0.16cm] 
    \nl  Compute temporary local coupling variables  \label{alg::hdqp::decoupled}
        \begin{equation}\label{eq::hdq::decoupled}
            \bar x_\ell =  \left(\overline{
            B}_\ell^k\right)^{-1}\ \left(\overline{B}^k_\ell x^k_\ell -\overline{g}^k_\ell\right).
        \end{equation}
        with $x^k_\ell=A_\ell\chi^k_\ell$, which essentially solves the decoupled subproblems for all $\ell \in \mathcal{R}$.\\
    \nl Compute weighted average  \label{alg::hdqp::averaging}
        \begin{equation}\label{eq::hdq::average}
            \bar z = \left(\sum_{\ell\in\mathcal R}E_\ell^\top \overline{B}_\ell^k E_\ell\right)^{-1} \sum_{\ell\in\mathcal R}E_\ell^\top \overline{B}_\ell^k \bar x_\ell.
        \end{equation}
    \nl Update primal and dual variables for all\label{alg::hdqp::update} $\ell\in\mathcal{R}$ by
    \begin{subequations}\label{eq::hdsqp::recovery}
        \begin{align}
            \lambda^{k+1}_\ell &= \overline{B}^k_\ell \left(\bar x_\ell -E_\ell  \bar z\right),\label{eq::hdq::dual}\\
            \chi_\ell^{k+1} &= \left(B_\ell^k\right)^{-1}\left(B_\ell^k \chi_\ell^k-g^k_\ell-A_\ell^\top \lambda_\ell^{k+1}\right).\label{eq::hdq::primal}
        \end{align} 
    \end{subequations}
}
\caption{Hypergraph-based Distributed Sequential Quadratic Programming (HDSQP)}\label{alg}
\end{algorithm}

Due to the positive definiteness of approximated Hessians based on the Levenberg-Marquardt method, the \acrshort{hdqp} can converge to the global minimizer of the condensed \acrshort{qp} subproblem~\eqref{eq::dsqp::subqp::reduced} in one iteration for saving computational and communication efforts. In the end, based on dual variable $\lambda$ provided by the inner algorithm \acrshort{hdqp}, the new full-dimensional iterate $\chi^{k+1}$ is updated by~\eqref{eq::hdq::primal}. Note that all the steps in Algorithm~\ref{alg} can be executed in parallel, except for the weighted averaging (Line~\ref{alg::hdqp::averaging}).

\subsection{Local Convergence Analysis}
\label{sec::analysis}
%\textcolor{red}{In the present paper, $z$ only denotes the common values of the coupling variables $x$. Hence, we focus on the convergence behaviors of the primal iterates $\chi$ and $x$ in this section. Moreover, the convergence of Algorithm~\ref{alg} by using the Levenberg-Marquardt method is equivalent to the Gauss-Newton method in this study since $\varepsilon=10^{-10}$ in~\eqref{eq::lm}. }
Without loss of generality, we assume that the flat start can provide a good initial guess (Remark~\ref{remark::flat}) such that the present paper focuses on the local convergence of Algorithm~\ref{alg}. Here, local means that the initial iterate is located in a small neighborhood of a local minimizer, within which the solution has physical meaning.

In the following, we first analyze the convergence of the inner algorithm \acrshort{hdqp} for solving condensed \acrshort{qp} subproblems~\eqref{eq::dsqp::subqp::reduced}, and prove that it can converge to a global minimizer in one step. Then, regardless of the inexactness caused by condensing subproblems, we prove that Algorithm~\ref{alg} can converge with a locally quadratic convergence rate when the Levenberg-Marquardt method is used to approximate Hessians.

%\acrshort{hdqp} is implemented to solve the condensed \acrshort{qp} subproblems~\eqref{eq::dsqp::subqp::reduced}. 
%By substituting \eqref{eq::hdq::average} into %\eqref{eq::hdq::dual}, we have
%\begin{equation}\label{eq::hdq::dual::matrix}
%    \lambda^{k+1}=\overline{B}^k M^k x
%\end{equation}
%with $M^k =I - E\left(E^\top \overline{B}^k E\right)^{-1}E^\top \overline{B}^k$.
%\begin{lem}[\cite{papastaikoudis_hypergraph_2022}]\label{lem::projection}
%    Matrix $M^k$ is a projection matrix, i.e.,
%    \begin{equation}
%        \left(M^k\right)^2 = M^k.
%    \end{equation}
%\end{lem}            % Lemma
%With this property, \acrshort{hdqp} has a fast convergence rate.
\begin{proposition}\label{thm::convergence::hdq}
Let the Levenberg-Marquardt method be used to evaluate $B_{\ell}^k$ such that the condensed \acrshort{qp} subproblem~\eqref{eq::dsqp::subqp::reduced} is strongly convex, then $$(x_\ell^{k+1}:=A_\ell\chi^{k+1},\; z^{k+1},z^{k+1}:=\bar z,\; \lambda_\ell^{k+1})$$  
given by Algorithm~\ref{alg} solves~\eqref{eq::dsqp::subqp::reduced} at iteration $k$.
\end{proposition}
The detailed proof appears in Appendix~\ref{app:prop}. As we know that the AC power flow equation~\eqref{eq::pf::equation} is sufficiently smooth, the objective $f$ is twice-Lipschitz continuously differentiable, i.e., 
there exists a constant $L>0$
\begin{equation}\label{eq::lipschitz}
    \frac{\norm{\nabla f(\chi)-\nabla f(\chi^{\ast})}}{\norm{\chi-\chi^{\ast}}}=\norm{\nabla^2 f(\tilde{\chi})} \leq L
\end{equation}
with $\tilde{\chi} = \chi - t(\chi-\chi^{\ast})$ for some $t\in(0,1)$. Moreover, since the optimal solution is feasible to~\eqref{eq::pf::equation}, we have zero-residual $r^*=0$ at the optimizer and 
\begin{equation}\label{eq::Q}
    Q^* = Q(\chi^\ast) =0.
\end{equation}
Locally, we can thus, have $\|Q(\chi)\| = \mathcal O(\|\chi-\chi^*\|)$. Before we establish the local convergence result of Algorithm~\ref{alg}, we introduce the definition of regular KKT point for~\eqref{eq::pf::ls::hypergraph}.
\begin{definition}[Regular KKT point of~\eqref{eq::pf::ls::hypergraph}]
A KKT point of~\eqref{eq::pf::ls::hypergraph} is called regular if \acrfull{sosc} and \acrfull{licq} hold at the KKT point~\cite{nocedal2006numerical}.
\end{definition}
\begin{remark}
    Due to~\eqref{eq::Q}, Hessian is equivalent to Gauss-Newton approximation at a local minimizer $\chi^*$. Moreover, the Jacobian matrix of the power flow equations~\eqref{eq::pf::equation} is always full-row rank in practice. Thus, \acrshort{sosc} is satisfied for the problem~\eqref{eq::pf::ls::hypergraph}.
\end{remark}
Additionally, the coupling introduced in Problem~\eqref{eq::pf::ls::hypergraph} is based on the hypergraph, we have \acrshort{licq} hold for the coupled affine equality constraints~\eqref{eq::coupling}, i.e., $E$ is full row rank. As a result, a \acrshort{kkt} point for the problem~\eqref{eq::pf::ls::hypergraph} is regular. 
\begin{theorem}
\label{thm::gn::rate}
    Let the minimizer $(\chi^{\ast},\lambda^{\ast})$ satisfy \acrshort{sosc} such that $(\chi^{\ast},\lambda^{\ast})$ is a regular \acrshort{kkt} point,
    let the parameter $\varepsilon$ be sufficiently small. Then, for solving the problem~\eqref{eq::pf::ls::hypergraph}, the iterates $\chi$ of Algorithm~\ref{alg} converges locally with a quadratic convergence rate.
\end{theorem}
If the exact Hessian is used in the \acrshort{qp} subproblems~\eqref{eq::subproblem}, the corresponding $p^\textsc{n}$ can be viewed as a standard Newton step of the original least-squares problem~\eqref{eq::pf::ls::original}. Similar to \cite{nocedal2006numerical}, we have
\begin{align}\label{eq::newton}
    %\norm{x^k+p^\textsc{n}-x^*} &\leq \hat{L}\norm{x^k-x^*}^2\\
    \norm{p^\textsc{n}} &\leq  \norm{\chi^k-\chi^*} + \norm{\chi^k+p^\textsc{n}-\chi^*}\notag\\
    &\leq  \norm{\chi^k-\chi^*} + \hat{L}\norm{\chi^k-\chi^*}^2,
\end{align}
where $\hat{L} = L \norm{(\nabla^2 f^*)^{-1}}$ and $L$ is the Lipschitz constant for $\nabla^2 f$ for $\chi$ near $\chi^*$. Since $\varepsilon$ is sufficiently small, the Levenberg-Marquardt method shares the same properties with the Gauss-Newton method. Hence, the analysis in the following is based on the Gauss-Newton method. The detailed proof is given in Appendix~\ref{app:rate}.

%$\omega_1 = \norm{\left(B^k\right)^{-1}}\cdot\norm{Q^k}$ and\\ $\omega_2 =\hat{L}\left(\norm{\left(B^k\right)^{-1}}\cdot\norm{Q^k}+1\right)$. 
%\begin{equation}
%    \omega_1 = \norm{\left(B^k\right)^{-1}}\cdot\norm{Q^k} \textrm{ and }\omega_2 =\hat{L}\left(\norm{\left(B^k\right)^{-1}}\cdot\norm{Q^k}+1\right)
%\end{equation}

\section{Numerical Case Study}\label{sec::simulation}

In this section, we illustrate the performance of the proposed distributed approach to solve AC \acrshort{pf} problems and compare it with the state-of-art \acrshort{aladin} algorithm.
\subsection{Implementation}
\begin{figure}[htbp!]
    \centering
    \includegraphics[width=\linewidth]{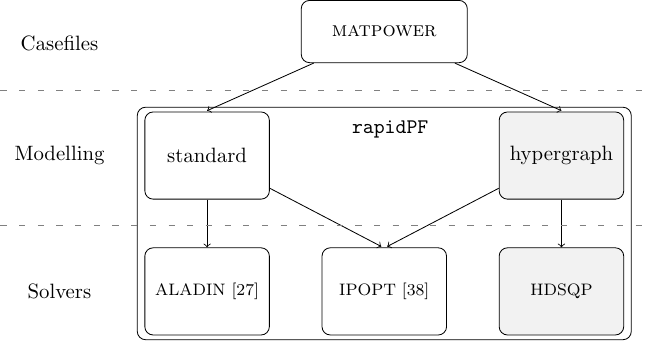}
    \caption{The open-source toolbox \acrshort{rapidpf}}
    \label{fig::rapidPF} 
\end{figure}

\change{The framework presented in this paper is implemented in \matlab-R2021a, and both the hypergraph-based problem and the proposed \acrshort{hdsqp} algorithm are provided in the \acrshort{rapidpf} toolbox~\footnote{Toolbox Documentation: \url{https://xinliang-dai.github.io/rapidPF/}}.}
%\footnote{Open-source toolkit: \url{https://github.com/xinliang-dai/rapidPF}}, a fully \matpower-compatible software that facilitates the laborious task on distributed AC \acrshort{pf} problems~\cite{muhlpfordt2021distributed,dai2022rapid}. 
As shown in Fig.~\ref{fig::rapidPF}, the toolbox allows users to combine multiple \matpower casefiles~\cite{zimmerman2010matpower} into a single merged casefile, formulate AC \acrlong{pf} problems as distributed optimization problems and then solve the problems by distributed approaches. Compared with previous work~\cite{muhlpfordt2021distributed,dai2022rapid}, the toolbox can reformulate the problems with a communication structure corresponding to a hypergraph
\eqref{eq::pf::ls::hypergraph} and solve them by the proposed approach~\acrshort{hdsqp}. Additionally, the \ipopt solver~\cite{wachter2006implementation} is used for calculating reference solutions, and the \casadi toolbox~\cite{andersson2019casadi} is used to compute exact Hessians for analysis.

The case studies are carried out on a standard desktop computer with \texttt{Intel\textsuperscript{\textregistered} %Core\texttrademark\, 
i5-6600K CPU @ 3.50GHz} and 16.0 \textsc{GB} installed \textsc{ram}. Following~\cite{dai2022rapid}, three benchmarks are created by using the \acrshort{rapidpf} toolbox based on IEEE standard test systems. The \acrshort{aladin} algorithm in the same toolbox is used for comparison. For a fair comparison, both \acrshort{aladin} and \acrshort{hdsqp} algorithms are initialized with a flat start. The computational time is estimated under the assumption that all subproblems are solved in parallel, and the time spent on exchanging sensitivities information is not taken into consideration.  %Case 1 merges a 30-bus system, a 14-bus system, and a 9-bus system into one grid, case 2 combines a 300-bus system and a 118-bus system, case 3 consists of 10 118-bus systems. }.

\begin{table*}[hbtp!]
    \centering
%\footnotesize
    \renewcommand{\arraystretch}{1.8}	
    \caption{Numerical Comparison: \acrshort{hdsqp} vs \acrshort{aladin}}\label{tb::results}
  \begin{tabular}{rrrrrrrcccc}  
  \hline
    Case & $n^\textrm{bus}$ &$n^\textrm{reg}$& $n^\textrm{state}$ & $n^\textrm{cpl}$ & Algorithm & Iterations& Time [$10^{-2}$s]& $\norm{\chi-\chi^*}$ & $\left|f(\chi)\right|$  & Primal Residual\\[0.1cm]\hline
    \multirow{2}{*}{1}&\multirow{2}{*}{53}    &  \multirow{2}{*}{3}    & \multirow{2}{*}{232}& \multirow{2}{*}{40}& \acrshort{aladin}&$4$&$2.33$&$6.13\times 10^{-09}$&$5.05\times 10^{-12}$&$3.61\times 10^{-9}$\\[0.1cm]
    &&&&&\acrshort{hdsqp}&$6$&$1.34$&$4.08\times10^{-08}$&$1.52\times10^{-20}$&$4.03\times 10^{-8}$\\[0.1cm]
    \multirow{2}{*}{2}&\multirow{2}{*}{418}    &  \multirow{2}{*}{2}   & \multirow{2}{*}{1684}& \multirow{2}{*}{24}& \acrshort{aladin}&$4$&$7.17$&$1.82\times10^{-09}$&$6.37\times10^{-09}$&$7.32\times 10^{-9}$\\[0.1cm]

    &&&&&\acrshort{hdsqp}&$10$&$5.68$&$1.39\times10^{-06}$&$3.04\times10^{-21}$&$1.37\times 10^{-9}$\\[0.1cm]
    \multirow{2}{*}{3}&\multirow{2}{*}{1180}   &  \multirow{2}{*}{10}  & \multirow{2}{*}{4764}& \multirow{2}{*}{88}& \acrshort{aladin}&$5$&$7.37$&$4.36\times 10^{-11}$& $2.59\times10^{-11}$&$2.56\times 10^{-9}$\\[0.1cm]
    &&&&&\acrshort{hdsqp}&$6$&$3.17$&$2.17\times10^{-08}$&$1.21\times10^{-20}$&$2.18\times 10^{-8}$\\[0.1cm]\hline
  \end{tabular}
\end{table*}
\begin{figure*}[htbp!]
    \centering
    \includegraphics[width=\textwidth]{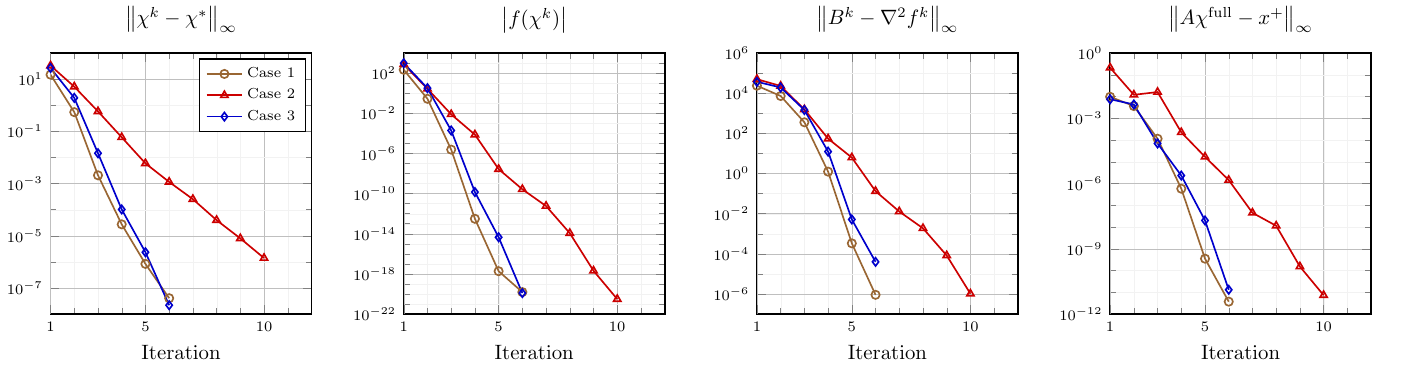}
    \caption{Convergence behavior of \acrshort{hdsqp}}\label{fig::convergence::hdsqp} 
\end{figure*}

\begin{remark}[Initialization]
\label{remark::flat}
   AC \acrshort{pf} problems are usually initialized with a flat start, where all voltage angles are set to zero, and all voltage magnitudes are set to 1.0 p.u.~\cite{frank2016introduction}. The initialization strategy has been demonstrated numerically that it can provide a good initial guess for the distributed approach in practice~\cite{muhlpfordt2021distributed,dai2022rapid}.
\end{remark}

\subsection{Case Studies}
Three test cases are studied in the present paper, as shown in Table~\ref{tb::results}. Note that $n^\textrm{reg}=|\mathcal R|$, $n^\textrm{state}$ and $n^\textrm{cpl}$ represent the number of regions, i.e., the cardinality of $\mathcal R$, the dimension of the state variables $\chi$ and the coupling variables $x$, respectively. To illustrate the convergence performance of the \acrshort{hdsqp} algorithm, we introduce four quantities, i.e.,
\begin{enumerate}
    \item deviation of state iterates to the minimizer $\norm{\chi^k-\chi^*}$,
    \item power flow residual $\left|f(\chi^k)\right|$,
    \item error of Levenberg-Marquardt Approximation
    $$\norm{B^k-\nabla^2f^k},$$
    \item inexactness caused by condensing \acrshort{qp} subproblems, i.e., deviation between the solution $x^+$ to condensed \acrshort{qp} subproblems~\eqref{eq::dsqp::subqp::reduced} with the solution $\chi^\textrm{full}$ to the full-dimensional \acrshort{qp} subproblems~\eqref{eq::subproblem}
    $$\|A\chi^k -x^k \|\quad \text{with}\;\;A=\text{diag}\{A_\ell\}_{\ell\in\mathcal R}.$$
\end{enumerate}

For a fair comparison, all the problems are initialized with a flat start. Table~\ref{tb::results} shows that the proposed \acrshort{hdsqp} can converge rapidly to a very highly accurate solution regarding the deviation of state variables, power flow residuals, and primal residual $x-E\,z$. Furthermore, the convergence behavior of the proposed algorithm for all three test cases is presented in Fig.~\ref{fig::convergence::hdsqp}. As the primal iterates $\chi^k$ approach, the minimizer $\chi^*$, Levenberg-Marquardt Approximation, and the solution to the \acrshort{qp} subproblems by using Schur decomposition become more accurate. Case 1 and Case 3 share almost the same performance and converge in 6 iterations, while Case 2 takes more iterations. The relative lower accuracy of the Schur complement possibly slows down the overall convergence rate of the proposed \acrshort{hdsqp} approach.

\subsection{HDSQP vs. ALADIN}

The \acrshort{aladin} algorithm used for comparison is a Gauss-Newton-based variant tailored to deal with AC power flow problem in~\cite{dai2022rapid}. As discussed in~\cite{dai2022rapid}, this \acrshort{aladin} variant has been illustrated that it outperforms the other existing state-of-art distributed approaches. Therefore, to illustrate the effectiveness of Algorithm~\ref{alg}, we compare it with this Gauss-Newton \acrshort{aladin} variant. 
\begin{figure}[htbp!]
    \centering
    \includegraphics[width=\linewidth]{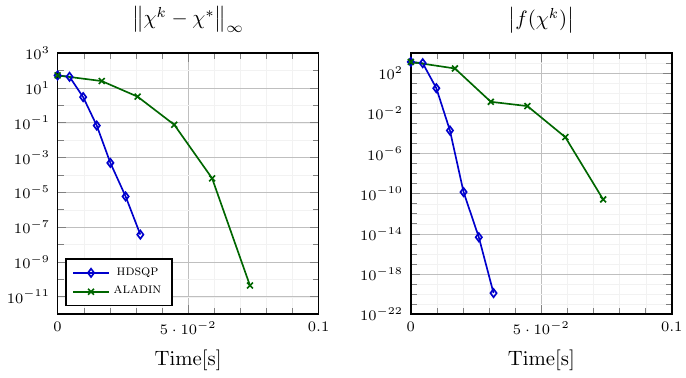}
    \caption{Comparison of different algorithms for Case 3}\label{fig::convergence::compare} 
\end{figure} 
First of all, let us have a close look at the computational complexity per iteration of both algorithms. The Gauss-Newton \acrshort{aladin} proposed in~\cite{dai2022rapid} requires $$\underbrace{\mathcal{O}(\sum_\ell( n^\textrm{state}_\ell)^3)}_{\text{parallelizable}}+\underbrace{\mathcal{O}(( n^\textrm{state})^3)}_{\text{consensus}}$$
float operations per iteration while \acrshort{hdsqp} needs 
$$\underbrace{\mathcal{O}(\sum_\ell( n^\textrm{state}_\ell)^3)}_{\text{parallelizable}}+\underbrace{\mathcal{O}(( n^\textrm{cpl})^3)}_{\text{consensus}}.$$
Here, one can see that the parallelizable computational complexity are same for both approaches, while the proposed Algorithm~\ref{alg} is much cheaper in the consensus part. 
This is because $n^\textrm{state}\gg n^\textrm{cpl}$ always holds for electric power systems in practice, as shown in Table~\ref{tb::results}.  Fig.~\ref{fig::convergence::compare} displays the convergence behaviors of both algorithms for Case 3. Although \acrshort{aladin} can converge with one iteration faster, the proposed \acrshort{hdsqp} has a much shorter computing time per iteration, benefitting from the fast convergence of the inner algorithm and thus surpasses \acrshort{aladin} in terms of total computing time. 

Table~\ref{tb::results} presents simulation results for comparison between the \acrshort{aladin} and the \acrshort{hdsqp} algorithms. For all three test cases, computational speedup of \acrshort{hdsqp} can be $30\%$ to $50\%$ compared with \acrshort{aladin} as presented in Table~\ref{tb::results}. Furthermore, \acrshort{hdsqp} requires only condensed Hessian in the centralized step, while \acrshort{aladin} requires both the first and the second order derivatives from all subproblems to solve a full-dimensional coupled \acrshort{qp} problem in the coordinator. This could further slow the total running time during parallel computing. Consequently, the proposed \acrshort{hdsqp} outperforms \acrshort{aladin} in solving AC \acrshort{pf} problems in aspects of computing time and communication effort.

\section{Conclusion and Outlook}\label{sec::conclusion}
The present paper proposes a distributed approach, \acrfull{hdsqp}, for solving power flow (PF) problems. \change{By introducing the hypergraph theory, the QP subproblems can be solved efficiently by the inner algorithm, i.e., the \acrfull{hdqp}~\cite{papastaikoudis_hypergraph_2022}.} A mathematical proof is provided that the inner algorithm \acrshort{hdqp} can converge in one iteration, and the local convergence rate of the proposed \acrshort{hdsqp} can achieve quadratic by implementing the Levenberg-Marquardt method to approximate Hessians. \change{Simulation results and analysis of the computational complexity demonstrate that the proposed algorithm outperforms the state-of-the-art distributed algorithm in terms of computing time for small- and medium-sized power grids at the cost of slightly increased iterations. Moreover, the numerical tests are added to the open-source toolbox \acrshort{rapidpf}.
%\footnote{Open-source toolkit: \url{https://github.com/xinliang-dai/rapidPF}}}.
} 

One drawback of the proposed approaches is associated with the inner algorithm~\acrshort{hdqp}. The inner algorithm employs a weighted averaging technique that utilizes the Hessian matrix, resulting in solutions converging to an equilibrium point near the exact optimizer. This hinders the scalability and numerical robustness of the proposed \acrshort{hdsqp}. To address these limitations, future work could focus on tuning the Levenberg-Marquardt method or alternating between \acrshort{hdsqp} and \acrshort{aladin} to enhance the scalability and numerical robustness.

\bibliographystyle{ieeetr}
\bibliography{mybib} 

%The present paper proposes a \acrfull{hdsqp} approach for solving \acrfull{pf} problems in a distributed fashion, where a \acrfull{hdqp} is implemented to solve \acrshort{qp} subproblems of \acrshort{hdsqp}. Furthermore, a mathematical proof is provided that the inner algorithm \acrshort{hdqp} can converge in one iteration, and the convergence rate of the proposed \acrshort{hdsqp} can achieve quadratic when approaching the minimizer if the Levenberg-Marquardt method is implemented to approximate Hessians. Simulation results validate the convergence performance of the proposed algorithm and show that it surpasses the state-of-art \acrshort{aladin} algorithm in terms of computing time for small- and medium-sized power grids. Future directions of the work include scaling up the algorithm to larger power systems and tuning parameters of the Levenberg-Marquardt method to improve the numerical robustness of the proposed algorithm.

\appendix
\subsection{Proof of Proposition~\ref{thm::convergence::hdq}}
\label{app:prop}
Based on~\cite{papastaikoudis_hypergraph_2022},
we can prove that the new iterate $(x^{k+1},z^{k+1},\lambda^{k+1})$ satisfies the \acrshort{kkt} condition~\eqref{eq::kkt}. %for the condensed \acrshort{qp} subproblem~\eqref{eq::dsqp::subqp::reduced}. 
By substituting \eqref{eq::hdq::average} into \eqref{eq::hdq::dual}, we have
\begin{equation}\label{eq::hdq::dual::matrix}
    \lambda^{k+1}=\overline{B}^k M^k \bar x
\end{equation}
with $M^k =I - E\left(E^\top \overline{B}^k E\right)^{-1}E^\top \overline{B}^k$. Consequently, we have
\begin{align}\label{eq::proof::dual}
    E^\top \lambda^{k+1}= E^\top \overline{B}^k M^k \bar x = 0.
\end{align}
This satisfies the dual feasibility~\eqref{eq::kkt::subqp::dual}. By substituting \eqref{eq::hdq::decoupled} into~\eqref{eq::hdq::dual::matrix}, we have
\begin{align}\label{eq::lam}
    \lambda^{k+1} = \overline{M}^k \,\overline{b}^k,
\end{align}
with $\overline{M}^k = \overline{B}^k M^k \left(\overline{B}^k\right)^{-1}$ and $\overline{b}^k= \overline{B}^k x^k-\overline{g}^k$. Accordingly, we can rewrite~\eqref{eq::hdq::primal} in a condensed form
\begin{align}
    x^{k+1} =& \left(\overline{B}^k\right)^{-1}\ \left(\overline{B}^k x^k -\overline{g}^k -\lambda^{k+1}\right)\label{eq::proof::reduced::primal}\\
    =&\left(\overline{B}^k\right)^{-1}\ \left(\overline{B}^k x^k -\overline{g}^k -\overline{M}^k\;\overline{b}^k\right)\notag\\
    =&  \left(\overline{B}^k\right)^{-1}\ \left(I- \overline{M}^k\right)\overline{b}^k\notag\\
    =& \left(I- M^k\right)\ \left(\overline{B}^k\right)^{-1} \overline{b}^k\notag\\
    =& E\left(E^\top \overline{B}^k E\right)^{-1}E^\top\; \overline{b}^k,
\end{align}
and its common value
\begin{align}
    z^{k+1} =& \left(E^\top \overline{B}^k E\right)^{-1} E^\top \overline{B}^k E\left(E^\top \overline{B}^k E\right)^{-1}E^\top\; \overline{b}^k\notag\\
    =& \left(E^\top \overline{B}^k E\right)^{-1}E^\top\; \overline{b}^k=\bar z.
\end{align}
Thereby, primal feasibility~\eqref{eq::kkt::subqp::primal} is satisfied by
\begin{align}\label{eq::proof::primal}
    x_\ell^{k+1} = E_\ell\,z^{k+1},\;\ell\in\mathcal R.
\end{align}
Moreover, the condition~\eqref{eq::kkt::subqp::1} is trivially satisfied due to~\eqref{eq::proof::reduced::primal}. We have thus established the Theorem~\ref{thm::convergence::hdq} by combining~\eqref{eq::proof::dual}~\eqref{eq::proof::reduced::primal}~\eqref{eq::proof::primal}.
\subsection{Proof of Theorem~\ref{thm::gn::rate}}
\label{app:rate}
%\small
The deviation between the Newton step $p^\textsc{n}$ and the Gauss-Newton step $p^\textsc{gn}$ can be written as
\begin{align}
    p^\textsc{n}-p^\textsc{gn} =& \left(B^k\right)^{-1} \left( B^k p^\textsc{n} + \nabla f^k\right)\notag\\
    =& \left(B^k\right)^{-1} \left(B^k-\nabla^2 f^k\right) p^\textsc{n}\notag\\
    =& - \left(B^k\right)^{-1} Q^k p^\textsc{n}
\end{align}
As a result, we obtain the following inequality
\begin{align*}
    \norm{\chi^k+p^\textsc{gn}-\chi^*}\leq \;& \norm{\chi^k+p^\textsc{n}-\chi^*} + \norm{p^\textsc{n}-p^\textsc{gn}}\notag\\
    \leq \;& \omega_1 \norm{\chi^k-\chi^*} + \omega_2 \norm{\chi^k-\chi^*}^2,
\end{align*}
where
\begin{equation}
    \omega_1 = s^k\cdot q^k \;\textrm{and}\; \omega_2 = \hat{L}(s^k\cdot q^k+1)
\end{equation}
with bounded $s^k = \norm{\left(B^k\right)^{-1}}$ and $q^k = \norm{Q^k} = \mathcal{O}(\|\chi^k-\chi^*\|)$. The locally quadratic convergence rate of Algorithm~\ref{alg} can be, thus, established~\cite{nocedal2006numerical}.

\subsection{Anonyms}\label{sec::anonyms}
\begin{table}[hbtp!]
    \centering
%\footnotesize
    \renewcommand{\arraystretch}{1.2}	
%    \caption{Numerical Comparison: \acrshort{hdsqp} vs \acrshort{aladin}}\label{tb::results}
  \begin{tabular}{ll}  
  \hline
    Abbr. & Description \\\hline
    \acrshort{admm} & \acrlong{admm} \\
    \acrshort{aladin} & Augmented Lagrangian based Alternating Direction Inexact \\
    & Newton method \\
    \acrshort{app} & \acrlong{app} \\
    \acrshort{hdqp} & \acrlong{hdqp} \\
    \acrshort{hdsqp} & \acrlong{hdsqp} \\
    \acrshort{kkt} & \acrlong{kkt} \\
    \acrshort{licq} & \acrlong{licq} \\
    \acrshort{nlp} & \acrlong{nlp} \\
    \acrshort{ocd} & \acrlong{ocd} \\
    \acrshort{pf} & \acrlong{pf} \\
    \acrshort{dsos} & \acrlong{dsos} \\
    \acrshort{tsos} & \acrlong{tsos} \\\hline
  \end{tabular}
\end{table}

%\printglossaries[type=main,style=long,nonumberlist]
%\printnoidxglossaries%[type=main,nonumberlist]
%\printnoidxglossaries

\end{document}